\documentclass{lms}[12pt]

\usepackage{amssymb,theorem,amsmath}

\title{A variant of the Hales-Jewett Theorem}

\author{Mathias Beiglb\"ock}

\classno{05D10 (primary), 22A15 (secondary)}

\extraline{The author thanks the Austrian Science Foundation 
FWF for its 
support through projects no.\
S8312 and no.\ S9612.}

\newtheorem{theorem}{Theorem}         
\newtheorem{corollary}{Corollary}
\newtheorem{lemma}{Lemma}

\newnumbered{assertion}{Assertion}    
\newnumbered{conjecture}{Conjecture}  
\newnumbered{definition}{Definition}
\newnumbered{hypothesis}{Hypothesis}
\newnumbered{remark}{Remark}
\newnumbered{note}{Note}
\newnumbered{observation}{Observation}
\newnumbered{problem}{Problem}
\newnumbered{question}{Question}
\newnumbered{algorithm}{Algorithm}
\newnumbered{example}{Example}
\newunnumbered{notation}{Notation}

\newcommand{\Z}{\mathbb{Z}}

\newcommand{\N}{\mathbb{N}}
\newcommand{\F}{\mathcal{F}}
\renewcommand{\P}{\mathcal{P}}
\newcommand{\Fin}{\P_f}

\newcommand{\rest}{\upharpoonright}
\newcommand{\dom}{{\it dom}}

\renewcommand{\Cup}{\bigcup}
\renewcommand{\Cap}{\bigcap}
\newcommand{\nhat }[1]{\{1,2,\ldots,#1\}}
\newcommand{\ohat }[1]{\{0,1,\ldots,#1\}}
\renewcommand{\subset}{\subseteq}
\renewcommand{\supset}{\supseteq}

\newcommand{\ovl}{\overline}

\begin{document}
\maketitle
\begin {abstract}

It was shown by V.\ Bergelson 
  that any set $B\subset \N$
 with positive 
upper multiplicative density contains nicely intertwined arithmetic and
geometric progressions:
For each $k\in \N$ there exist $a,b,d\in \N$ such that
$ \left\{b(a+id)^j:i,j \in\nhat k\right\}\subset B. $
In particular one cell of each finite partition of $\N$ contains such
 configurations. 
We prove a Hales-Jewett type extension of this partition theorem.
\end{abstract}

\section{Introduction}

Van der Waerden's Theorem (\cite{refW}) states that for any finite coloring of 
$\N$ one can find arbitrarily long monochromatic arithmetic progressions.
In 1963 A.\ Hales and R.\ Jewett (\cite{refHJ}) gave a powerful abstract extension of 
van der Waerden's Theorem.

We introduce some notations to state their result. An 
\emph{alphabet} $\Sigma$
 is a finite nonempty set. A \emph{located word} $\alpha$ is a 
function from a finite set $\dom\, \alpha \subset \N$ to $\Sigma$.
The set of all located words will be denoted by $L(\Sigma)$.
Note that for located words $\alpha,\beta$ 
satisfying  $\dom \, \alpha \cap \dom \, \beta = \emptyset$, 
$\alpha \cup \beta $ is also located word. (Here it is convenient 
to view functions as sets of ordered pairs.)
By $\Fin(\N)$ we denote the set of all finite nonempty subsets
of $\N$. 
\begin{theorem}[Hales-Jewett]\label{hj}
Let $L(\Sigma)$ be finitely coloured. There exist 
$\alpha\in L(\Sigma)$ and $\gamma\in \Fin(\N)$ 
 such that 
$\dom\, \alpha \cap \gamma = \emptyset$ and
$ \{\alpha \cup \gamma\times\{s\}:s\in\Sigma\} $ is monochrome.
\end{theorem}
The term
 $\alpha \cup \gamma \times \{s\}$ may be viewed as an analogue of 
the expression $a+c\cdot s$. In particular we use  
$\alpha \cup \gamma \times \{s\}$ for what should
rigorously be $\alpha \cup (\gamma \times \{s\})$. Configurations of the form $ \{\alpha \cup \gamma\times\{s\}:s\in\Sigma\} $ are often called \emph{combinatorial lines}.

We will explain shortly how van der Waerden's Theorem can be
derived from the Hales-Jewett Theorem. Let $\Sigma =\ohat k$ and 
assume that $ \N$ is finitely coloured. Consider
 the map $f:L(\Sigma )\rightarrow \N, \alpha\mapsto 
1+ \sum_{t\in \dom\, \alpha} \alpha(t)$ and colour each $\alpha\in
L(\Sigma)$ with the colour of $f(\alpha)$. Pick $\alpha$ 
and $\gamma$
 according
to Theorem \ref{hj}. Let $a = f(\alpha)$ and $d=|\gamma|$.
Then for all $i \in \ohat k$, $a+id =f (\alpha\cup 
\gamma\times \{i\})$ and thus the arithmetic progression 
$\{a, a+d, \ldots, a+kd\}$ is monochrome.
\\

Hales-Jewett type extensions of various other Ramsey theoretic results 
have been obtained. We mention two very deep theorems in this 
style:
H.\ Furstenberg and Y.\ Katznelson (\cite{refFK})
gave a density version of the Hales-Jewett Theorem, which generalizes 
Szemer\'edi's Theorem (\cite{refS}). V.\ Bergelson and A.\ Leibman 
(\cite{refBL2}) proved
a polynomial Hales-Jewett theorem 
which extends the polynomial van
der Waerden Theorem.

As shown in \cite{refB}, Theorem 3.15, for every set $B\subset\N$ 
with positive upper multiplicative density (actually positive upper
multiplictive Banach density is enough) and each
$k\in\N$ there exist $a,b,d\in \N$ such that 
$\left\{b(a+id)^j: i,j\in\ohat k\right\}\subset B$. 
A direct consequence 
is the following partition theorem:
\begin{theorem}\label{unsers}
For any finite coloring of $\N$ and  $k\in\N$
there exist $a,b,d\in \N$ such that
$\left\{b(a+id)^j:i,j\in\nhat k\right\}$ is monochrome.
\end {theorem}
(See \cite{refBBHS}) for an algebraic proof of this result.)

The main theorem of this paper
 is an extension of the Hales-Jewett Theorem which 
is strong enough to yield Theorem \ref{unsers} but also implies some 
other corollaries.
(Call a family $\F$ of subsets of $\N$ \emph{partition regular} if
for any finite coloring of $\N$ there exists some $F\in \F$ which is 
monochrome.)

\begin{theorem}\label{maintheorem}
Let $\F$ be a partition
regular family
of finite  subsets of $\N$
which contains no singletons and let $\Sigma $ be a finite alphabet.
For any finite colouring of $L(\Sigma)$
there exist $\alpha\in L(\Sigma), \gamma\in \Fin (\N)$ and $F\in \F$
such that $\dom \, \alpha, \gamma$ and $ F$ are pairwise disjoint
and
$$\big\{\alpha \cup (\gamma \cup \{t\})
\times \{s\} : s\in \Sigma, t\in F \big\}$$
is monochrome.
\end{theorem}

 Similarly as the Hales-Jewett Theorem implies van der
Waerden's Theorem, Theorem \ref{maintheorem} can be
 applied to derive Theorem \ref{unsers}.
Assume that $\N$ is finitely coloured. Fix $k\in\N$, let
$\F=\big\{\{a,a+d,\ldots,a+kd\}:a,d\in \N\big\}$ be the set of
all $(k+1)$-term arithmetic progressions, put $\Sigma=\ohat k$ and
define $f:L(\Sigma)\rightarrow \N$ by $f(\alpha)=\prod_{t\in\dom\, 
  \alpha} t^{\alpha(t)}.$ Colour each $\alpha \in L(\Sigma)$ with the 
colour of $f(\alpha)$ and choose $\alpha, \gamma$ and $ F=\{a, a+d,
  \ldots, a+kd\}$ according to Theorem \ref {maintheorem}. 
Then for all $i,j \in \ohat k$, 
$$\textstyle{f\big(\alpha\cup(\gamma\cup \{a+id\})
\times\{j\}\big)=\prod_{t\in\dom\, \alpha}t^{\alpha(t)}
\cdot \prod_{t\in\gamma}t^j \cdot (a+id)^j =}$$
$$\textstyle{=\underbrace{
   \textstyle{\prod_{t\in\dom\,\alpha}t^{\alpha(t)}}}_{=\ovl b}
    \cdot\Big(\underbrace{ 
   \textstyle{ a\cdot \prod_{t\in\gamma}t }}_{=\ovl a} +
          i\cdot \underbrace{
  \textstyle {d\cdot\prod_{t\in\gamma} t}}_{=\ovl d} 
   \Big)^j}=\ovl b(\ovl a+i\ovl d)^j
  $$
has the same colour.

Note that we may replace $\F$ by an arbitrary partition regular
 family of finite subsets of $\N$ which contains no singletons.
In this way we see that there exist $a,r\in\N$ and $ F\in \F$ 
such that $\left\{ b(rt)^j:t\in F, j\in \ohat k\right\}$ is monochrome.
(This stronger version also follows  from the algebraic proof of 
Theorem \ref{unsers}, see \cite{refBBHS}, Corollary 4.3.) 

Another application of Theorem \ref{maintheorem} is that 
it allows to lower the assumption on the set which is coloured
in Theorem \ref{unsers}.

\begin{corollary}\label{appli}
Let $k,r\in\N$. There exists $K\in\N$ such that for all
 $A, D \in \N$ the following holds:  Whenever 
$$ S_K(A,D)=\big\{ (A+i_1D)(A+i_2D) \ldots (A+i_mD):
   m, i_1,i_2, \ldots ,i_m \in \ohat K\big\}$$
is partitioned into sets $B_1,B_2,\ldots, B_r$
one can find $b,a,d\in\N$ and $s\in\nhat r$ such that 
$\big\{b(a+id)^j:i,j\in \nhat k\big\}\subset B_s$.
\end{corollary}
\begin{proof} Let $\Sigma=\ohat k$ and $\F=\big\{\{a,a+d,
\ldots,a+kd\}:a,d\in\N\big\}$.
 Using a standard compactness argument (cf. \cite{refHS}, section 5.5) 
we see that for fixed $r\in \N$
 there exists some $N\in\N$ such that for any
colouring of $L(\Sigma)$ into $r$ colours one can 
 choose $\alpha, \gamma$ and 
$F=\{a,a+d, \ldots ,a+kd\}$ according to Theorem \ref{maintheorem}
and additionally require $\dom\, \alpha, \gamma, F\subset \nhat N$.
Let $f:L(\Sigma)\rightarrow \N, f(\alpha)=\prod_
{t\in \dom \, \alpha} (A+tD)^{\alpha(t)}$ and put 
$K=kN.$ Thus $f(\alpha)\in S_K(A,D)$
 for all $\alpha \in L(\Sigma)$,
$\dom\, \alpha\subset \nhat N$. 
For each $s\in \nhat r$ 
let $C_s$ be the set of 
all  $\alpha\in L(\Sigma), \dom\, \alpha \subset\nhat N$
such that $f(\alpha)\in B_s$. Choose $s\in \nhat r$ and
 $\alpha,\gamma, F=\{a,a+d,\ldots,a+kd\}$
such that $ \dom\, \alpha,\gamma,F$ are 
pairwise disjoint subsets of
$\nhat N$ and 
$\alpha \cup (\gamma \cup \{a+id\})\times\{j\}\subset C_s$ 
 for all choices of $i,j\in \ohat k$.
Applying $f$ it follows that 
$\big\{\ovl b(\ovl a + i \ovl d)^j:i,j\in\ohat k\big\}\subset B_s$, where 
$\ovl b = \prod_{t\in\dom \, \alpha} (A+tD)^{\alpha(t)},
\ovl a = (A+aD) \prod_{t\in\gamma} (A+tD),$
and 
$\ovl d =d D \prod_{t\in \gamma} (A+tD).$
\end{proof}

Many structures considered
 in Ramsey theory have the nice property to be `unbreakable'
in the sense that if a sufficiently large structure of a certain 
type is partitioned into a specified in advance number of cells, 
at least one cell will again contain a large strucure of the same type.
For example it follows from van der Waerden's Theorem that for all 
$r,k\in\N$
there exists some $K\in\N$, such that whenever a $K$-term arithmetic progression
is partitioned into $r$ cells at least one cell contains
a $ k$-term arithmetic progression.
Unfortunately the intertwined additive-multiplicative
 structures in Theorem
\ref{unsers} do not posses this property as shown in \cite[Theorem 4.9]{refBBHS}
that for all $K\in\N$ 
there exist $ A, D\in\N$ 
and a partition $B_1\cup B_2=\left\{(A+iD)^j :i,j\in\ohat K\right\}$ 
such that no $B_i$ contains a configuration of the
form $\left\{ b(a+id)^j:i\in \{0,1,2\}, j\in \{0,1\}\right\}$. 
Corollary \ref{appli} gives a vague hint that 
the situation could be different
for sets 
$$S_k(a,b,d)=\big\{b(a+i_1d)(a+i_2d)\ldots(a+i_m d):m,i_1, i_2,\ldots,
   i_m \in\ohat k\big\},$$
where $a,b,d\in\N$. 

\begin{question}\label{Q1}
Fix $k,r\in \N$.
Does there exist some $K\in \N$ such that 
for all $A,B,D\in\N$ and any partition $\Cup_{s=1}^r B_s
  = S_K(A,B,D)$ one 
can find $a,b,d,k\in\N$ and $s\in\nhat r$ such that $S_k(a,b,d)\subset B_s$?
\end{question}

Note added in revision: Imre Leader (private communication) has found an example answering  Question \ref{Q1} in the negative. 
\section{Preliminaries}\label{pre}

We give a short outline of the ideas behind our proofs.
Lemma \ref{piecewiselemma} states that any large enough, that is, piecewise syndetic set contains structures as in Theorem \ref{maintheorem}. The argument behind this is abstract but simple in nature.
Then one shows that for any colouring of $\Sigma(L)$ there is a piecewise syndetic set $A$ of combinatorial lines having the same colour. Thus Lemma \ref{piecewiselemma} can be applied to $A$ to yield Theorem \ref{maintheorem}.

It was shown by Furstenberg and Glasner (\cite{refFG}) that any piecewise syndetic set of integers contains a piecewise syndetic set of arithmetic progressions. Their proof is based on the theory of compact semigroups applied to the Stone-\v Cech compactification $\beta \Z$ (in the form of the enveloping semigroup of an appropriate dynamical system)  and so far no purely elementary proof of this strong version of van der Waerden's Theorem is known. In \cite{refBH} similar statements are proved for various notions of largeness in an abstract semigroup  $S$. This extension is possible since - to some extend - the algebraic theory of $\beta S$ resembles  that of $\beta \Z$. In our case we have to work with the partial semigroup  of  located words  $L(\Sigma)$ respectively with its Stone-\v Cech compactification. We reassure the reader that under the right technical precautions this setting is closely related to the one encountered when facing  $\beta \Z$ (which of course is an esoteric enough object).
Below we  review the most important facts needed for  the proof of Theorem \ref{maintheorem} in Section \ref{pro}.

A \emph{partial semigroup}  
 $(S,\cdot)$ (as introduced in \cite{refBBH}) is a  set $S$
together with a binary operation $\cdot$ that maps a 
subset of $S\times S $ into $ S$ and satisfies the 
associative law $ x \cdot (y \cdot z)=(x\cdot y) \cdot z$
in the sense that if either side is defined, so is the other
and they are equal. (The notion of a partial semigroup was introduced in 
\cite{refBBH}.)

Given a partial semigroup $(S,\cdot)$ and $x\in S$
 let $\phi(x)=\{y\in S: x\cdot y\ \mbox{is defined}\}$.
$(S,\cdot)$ is called \emph{adequate} if
$\Cap_{x\in F}\phi(x)\neq\emptyset$  for all
finite nonempty $F\subset S$. 

We shall deal mostly with the adequate semigroup
$(L(\Sigma),\uplus)$ of located words over an
alphabet $\Sigma$ where we let 
$\alpha \uplus \beta=  \alpha\cup \beta $
for $\alpha,\beta\in L(\Sigma)$ if $\dom\, \alpha
\cap \dom \, \beta = \emptyset$ and leave $\alpha \uplus \beta$
 undefined otherwise.

Given a discrete space $S$ we take the 
Stone-\v Cech compactification $\beta S$ of $S$ to be the 
set of all ultrafilters on $S$, the points of $S$ being 
identified with the principal ultrafilters.
Given a set $A\subset S$ let $ \ovl A =\{p\in \beta S: A\in p\}$.
The family $\{\ovl A: A\subset S\}$ forms a clopen basis of 
$\beta S$. A  semigroup structure on $S$  can be extended to
$\beta S$. The resulting algebraic properties of 
$\beta S$ turn out to be extremly useful in Ramsey theory, see \cite{refHS}
for an extensive treatment  of this topic and related material.

In the case of partial semigroups a slightly different approach 
is used. For an adequate partial semigroup $(S,\cdot)$
let $\delta S=\Cap_{x\in S} \ovl{\phi (x)}\subset \beta S$.
For $x\in S$ and $A\subset S$ we let
$x^{-1}A=\{y\in \phi (x): x\cdot y\in A\}$. (Note that 
this corresponds to the usual definition if $(S,\cdot )$
is a group.)
For $p,q\in \delta S$ put
$$ p\cdot q = \{A\subset S:\{x\in S:x^{-1}A \in q\} \in p\}.$$ 
Then $(\delta S,\cdot)$ is a \emph{compact right topological
semigroup}, i.e.\ $\delta S$ is
compact with the topology inherited from
$\beta S$ and for each $q\in \delta S$ the map 
$p\mapsto p\cdot q$ is continuous (\cite{refBBH}, Proposition 2.6).

We will need some facts on
the algebraic properties of compact right topological
semigroups: Every compact right topological semigroup
$(T,\cdot)$ posses an 
\emph{idempotent}, i.e.\ there is some $e\in T$ such that 
$e\cdot e = e$. The set of all idempotents of $T$ 
is partially ordered if we let $e\leq f \ \Longleftrightarrow 
\ e\cdot f = f\cdot e =e$. Idempotents which are 
minimal with respect to this 
ordering are called \emph{minimal idempotents}.  
For any idempotent $f\in T$ there exists some
minimal idempotent $e\in T$ such that 
$e\leq f$. A different  characterisation
of the minimal idempotents in $T$ can be given
via the ideals of $T$. (A nonempty set $I\subset T$ is called an ideal
if $I\cdot T \cup T\cdot I\subset I.$) 
Every compact right topological
semigroup $T$ has a smallest ideal $K(T)$. An
idempotent $e\in T$ is a minimal
idempotent if and only if $e\in K(T)$. (See \cite{refHS} for a 
comprehensive  introduction to the theory of compact right 
topological semigroups.)

\section{Proof of the main Theorem}\label{pro}

Let $(S,\cdot)$ be a partial semigroup, and let  $\F$ be
 a family of subsets 
of $S$. $\F$ is called invariant if for all $s\in S$ and $F\in\F$
the following conditions hold:
\begin{enumerate}
\item If $s\cdot F =\{ s\cdot f: f\in F\}$ is defined, then $s\cdot F \in \F$.
\item If $F\cdot s=\{f\cdot s :f\in F\}$ is defined, then $F\cdot s \in \F$.
\end{enumerate}

$\F$ is \emph{adequately partition regular} if for any 
finite set $G\subset S$ and  any finite partition of 
$\Cap_{x\in G}\phi(x)$ 
there
exists some $F\in\F$ which is entirely contained in one cell of
the partition. 

If $(S,\cdot)$ is a partial semigroup,
 $A\subset S$ is called \emph{piecewise syndetic} if there
exists some $p\in K(\delta S)$ such that $A\in p$. 
Piecewise syndetic subsets of semigroups admit a simple combinatorial characterisation (see \cite[Theorem 4.40]{refHS}). For instance $A\subset (\Z,+)$ is piecewise syndetic if and only if $\Cup_{n=1}^r A-n$ contains arbitrarily long intervals.
However the combinatorial and the algebraic definition of piecewise syndetic lead to different concepts in the case of partial semigroups. (See \cite{refM,refM2} for a detailed analysis of this phenomenon.) For our intended application in the proof of Theorem \ref{maintheorem} the algebraic version is appropriate.  

It is a fairly easy combinatorial exercise to show that van der Waerden's Theorem corresponds to the fact that piecewise syndetic sets in $\Z$ contain arbitrarily long arithmetic progressions. The following lemma translates this idea to our setting.
\begin{lemma}\label{piecewiselemma}
Let $S$ be an adequate partial semigroup,
let $\F$ be an adequately partition regular invariant family of finite
subsets of $S$ and assume that $A\subset S$ is piecewise syndetic.
Then there exists $F\in \F$ such that $F\subset A$.
\end{lemma}
\begin{proof}
Let $P=\{p\in \delta S: \mbox{For all } A\in p \mbox{ exists } F \in \F
 \mbox{ such that } F\subset A\}$.
It is sufficient to show that $P$ is an ideal of $\delta S$, since 
this implies $ P\supset K(\delta S)$. 

Let $G$ be a finite subset of $ S$.  Since for any finite partition of 
$\Cap_{x\in G} \phi (x)$ there exists some $F\in \F$ 
which is entirely contained in one cell of the partition
we have that $P_G=\{p\in 
\ovl{\Cap_{x\in G} \phi (x)}: \mbox{For all } A\in p \mbox{ exists } F \in \F
 \mbox{ such that } F\subset A\}\neq \emptyset$ by \cite{refHS}, Theorem 
3.11. 

Each $P_G$ is closed, so by compactness of $\delta S,$ 
$P=\Cap_{G\subset S, |G|<\infty } P_G \neq \emptyset.$

To see that $P$ is a left ideal, let $p\in P$ and $q\in \delta S$.
Assume that $ A\in q\cdot p$, i.e.\ $ \{s:s^{-1}A\in p\}\in q.$  
Thus we may pick $s\in S$ such that $s^{-1}A\in p$. Since $p\in P$
there exists $F\in\F$ such that $F\subset s^{-1}A.$ This is equivalent to
$s F \subset A$ and $sF\in\F$ by the invariance of
 $\F$. 
Since $A$ was arbitrary we see that $q\cdot p\in P.$

To show that $ P$ is a right ideal, pick $p\in P, q\in\delta S$
and $ A\in p\cdot q$. Thus $\{s\in S :s^{-1}A\in q\}\in p$, so 
by the choice of $p$ we may pick $F\in \F$ such that $ 
F\subset \{s\in S:s^{-1}A\in q\}$. Since $F$ is finite
$\Cap_{s\in F}s^{-1} A\in q$, so pick $t\in \Cap_{s\in F}s^{-1} A$.
Then $F\cdot t\subset A$ and, once again by the invariance of
$\F$, $F\cdot t\in \F$.   Thus $p\cdot q \in P$.
\end{proof} 

 Let $v$ be a `variable' not in $\Sigma$.
$L(\Sigma\cup\{v \})$ is the set of all located words over
the alphabet $\Sigma\cup \{v\}$. 
 By $L(\Sigma ; v)= L(\Sigma\cup v) \setminus L ( \Sigma)$
 we denote the set  
of all located words over the alphabet $\Sigma \cup \{ v\}$ in 
which  $v$ occurs. The elements of $ L(\Sigma;\{v\})$ are often
called \emph{variable words}
 while one refers to elements of $L(\Sigma)$
as \emph{constant words}.
$L(\Sigma;v)$ is an ideal of $L(\Sigma\cup \{ v\})$  and
consequently $ \delta L(\Sigma;v )$ is an ideal of 
$\delta L(\Sigma\cup\{v\})$.

For $s\in \Sigma$ and $\alpha\in L(\Sigma\cup\{v\})$,
let $\theta_s(\alpha)$ be the result of replacing
each occurence of $v$ in $\alpha$ by $s$.
More formally, $\dom\, \theta_s(\alpha) =
\dom\, \alpha$ and for $t\in \dom\, \theta_s(\alpha)$
$$\theta_s(\alpha)(t)=\left\{\begin{array}{cl}
s & \mbox{if } \alpha(t)=v\\
\alpha(t) & \mbox{if } \alpha(t)\in\Sigma.
\end{array}\right.$$
The mapping $\theta_s :L(\Sigma\cup\{v\})\rightarrow L(\Sigma)$ 
gives rise to the  continuous extension $\tilde\theta_s:\beta 
L(\Sigma\cup\{v\})\rightarrow \beta L(\Sigma)$. 
Whenever $\alpha \uplus \beta $ is defined for 
$\alpha,\beta\in L(\Sigma)$, so is $\theta_s(\alpha)
\uplus \theta_s(\beta)$ and it equals $ \theta_s(\alpha\uplus \beta)$.
So $\theta_s:L(\Sigma\cup \{v\}) \rightarrow L(\Sigma)$ is
a \emph{homomorphism of partial semigroups}.
 Moreover $ \theta_s$ is surjective. By 
\cite{refBBH},
Proposition 2, this implies that the restriction
of $ \tilde \theta_s $ to $\delta L(\Sigma\cup\{v\})$ is a continuous 
homomorphism $ \delta L(\Sigma\cup\{v\})\rightarrow \delta L(\Sigma)$.
Since $ \theta_s\rest\ L(\Sigma)$ is the identity
the same holds true for $\tilde \theta_s \rest\ \delta L(\Sigma)$.

We are now able to give the proof of our main Theorem:
\begin{proof}[of Theorem \ref{maintheorem}]
Pick a minimal idempotent $e\in \delta L(\Sigma)$. 
Let
$M=\{ p\in\delta L(\Sigma\cup\{v\}):\tilde\theta_s(p)=e \mbox{ for all } 
  s \in \Sigma\}.$
Since, for each $s\in\Sigma$,
 $\tilde\theta_s$ is the identity on $\delta L(\Sigma)$ 
we have $\tilde\theta_s (e)=e.$ Thus $M$ is nonempty.
$M$ is the intersection of homomorphic
 preimages of the closed semigroup $
\{e\}$, thus it is a closed semigroup as well. 
Fix a minimal idempotent $q\in M$. We want to show
that $q$ is also minimal in $\delta L(\Sigma\cup \{v\})$. 
Pick a minimal idempotent $f\in \delta L(\Sigma\cup \{v\})$ 
 such that $f\leq q,$ i.e.\ $f\uplus q= q\uplus f =f$. 
Let $s\in \Sigma$. We have $
\tilde \theta_s(f)=\tilde\theta_s (f\uplus q )=
\tilde\theta_s(f) \uplus \tilde\theta_s(q)= 
\tilde\theta_s(f)\uplus e$ and, analoguously, 
$\tilde\theta_s(f)=e\uplus\tilde\theta_s(f)$, thus
$\tilde\theta_s(f)\leq e$. Since $e $ was chosen to be minimal   
in $L(\Sigma)$ it follows that $\tilde\theta_s(f)=e$. $s\in\Sigma$
was arbitrary, thus we have $f\in M$. Since $q$ is minimal in $M$ it 
follows that $q=f$, so $q$ is in fact a minimal idempotent in
$\delta L(\Sigma\cup\{v\})$.

By the ultrafilter property of $e$, choose a
 monochrome set $B\subset L(\Sigma)$ such that $B\in e$. 
Thus $\overline B$ is a neighbourhood of $e  $, so pick, 
by continuity of $\tilde\theta_s, s\in\Sigma$,
 a neighbourhood $\overline A$ of $q$ such 
that  
$\tilde \theta_s[\ovl A]\subset \overline B$
for all $s\in\Sigma$.
Then 
$A\in q\in\delta L(\Sigma\cup \{v\})$,
 so $A$ is piecewise syndetic in $L(\Sigma\cup\{v\})$ and
 $\theta_s[A]\subset B$ for all $s\in\Sigma$.

$\F$ is a partition regular family which contains no singletons. 
This implies that for any $m\in\N$ and  any finite colouring of 
$\{n\in \N:n> m \}$ there exists a monochrome set $F\in\F$. 
(Extend the colouring of $\{n\in\N:n> m\}$ to a colouring of 
$\N$ by giving all elements of $\nhat m$ new and mutually different 
colours. Any $F\in\F$ which is monochrome with respect to this colouring
is contained in $\{n\in\N:n>m\}$ since $F$ has more then one element.)

Consequently,
 $ \F'=\{\{\{(t,v)\}:t\in F\}: F\in \F\}$ is an
  adequately partition regular family 
of subsets of $L(\Sigma\cup \{v\})$ and thus
 $$\F''=\{\{\beta \uplus \{(t,v)\} : t\in F\}: \beta\in
    L(\Sigma;v), F\in \F, \dom\, \beta \cap F =\emptyset \} $$
 is an invariant
adequately partition regular 
family.  
So we may apply Lemma \ref{piecewiselemma} and pick
 $G\in \F''$ such that $ G\subset A$.
Each variable word $\beta\in L(\Sigma;v)$ can be written
in the form $\alpha\cup \gamma\times\{v\}$ for uniquely
determined
$\alpha\in L(\Sigma), \gamma\in \Fin(\N), 
\dom\, \alpha \cap \gamma=\emptyset.$
Hence $G$ is of the form 
$\{\alpha \uplus (\gamma\cup \{t\}) \times \{v\}:t\in F\}$, where
$\alpha \in L(\Sigma), \gamma\in\Fin(\N)$ and $F\in \F$.
It follows that for each $s\in \Sigma$ and each 
$t\in F,$ 
$$\theta_s(\alpha\uplus (\gamma\cup\{t\} )\times\{v\})=
  \alpha \cup (\gamma\cup\{t\})\times \{s\}\in B.$$
\end{proof}

It would be preferable to show Theorem \ref{maintheorem} in a purely elementary way. While it is not difficult to substitute Lemma \ref{piecewiselemma} by a combinatorial argument, the author does not know how to do this in the case of the other part of the proof.

\affiliationone{
   M.\ Beiglb\"ock\\
      TU Vienna\\
   Wiedner Hauptstr.\ 8-10\\
   1040 Vienna\\
   Austria\\
   \email{mathias.beiglboeck@tuwien.ac.at}}

\end{document}